\documentclass{article}
\usepackage{graphicx}
\usepackage{amscd}
\usepackage[matrix,arrow,curve]{xy}
\usepackage{amssymb,amsmath}
\usepackage{amsmath,amsfonts,amssymb,latexsym}
\newtheorem{theo}{Theorem}[section]
\newtheorem{mydf}{Definition}[section]

\newtheorem{cor}{Corollary}[section]
\newtheorem{lem}{Lemma}[section]

\newtheorem{exam}{Example}[section]
\def\P{\bf Proof. \rm}

\newcommand{\proofterminator}{\hfill\ensuremath{\square}}
\begin{document}
\title{Cofibrations in the Category of Fr\"olicher Spaces: Part I}
\author{{\small \underline{Brett Dugmore}} \\ {\small  Cadiz Financial Strategists (Pty) Ltd, Cape Town, South
Africa} \\ {\small Email: Brett.Dugmore@cadiz.co.za}\\ {\small \underline{Patrice Pungu Ntumba}}\\
{\small Department of Mathematics and Applied
Mathematics}\\{\small University of Pretoria}\\ {\small Hatfield
0002, Republic of South Africa}\\{\small Email:
patrice.ntumba@up.ac.za}\\}

\date{}
\maketitle

\begin{abstract} Cofibrations are defined in the category of
Fr\"olicher spaces by weakening the analog of the classical
definition to enable smooth homotopy extensions to be more easily
constructed, using flattened unit intervals. We later relate
smooth cofibrations to smooth neighborhood deformation retracts.
The notion of smooth neighborhood deformation retract gives rise
to an analogous result that a closed Fr\"olicher subspace $A$ of
the Fr\"olicher space $X$ is a smooth neighborhood deformation
retract of $X$ if and only if the inclusion $i: A\hookrightarrow
X$ comes from a certain subclass of cofibrations. As an
application we construct the right Puppe sequence.
\end{abstract}

\section{Preliminaries} The purpose of this section is to survey
brielfy the notion of Fr\"olicher spaces. Fr\"olicher spaces arise
naturally in physics, and do generalize the concept of smooth
manifolds. A \textsl{Fr\"olicher space}, or smooth space as
initially called by Fr\"olicher and Kriegl \cite{fk}, is a triple
$(X, \mathcal{C}_X, \mathcal{F}_X)$ consisting of a set $X$, and
subsets $\mathcal{C}_X\subseteq X^\mathbb{R}$,
$\mathcal{F}_X\subseteq \mathbb{R}^X$ such that
\begin{itemize}\item $\mathcal{F}_X\circ \mathcal{C}_X= \{f\circ
c|\ f\in \mathcal{F}_X, \ c\in \mathcal{C}_X\}\subseteq
C^\infty(\mathbb{R})$ \item $\Phi\mathcal{C}_X:= \{f: X\rightarrow
\mathbb{R}|\ f\circ c\in C^\infty(\mathbb{R})\ \mbox{for all $c\in
\mathcal{C}_X$}\}= \mathcal{F}_X$ \item $\Gamma\mathcal{F}_X:=
\{c: \mathbb{R}\rightarrow X|\ f\circ c\in C^\infty(\mathbb{R})
\mbox{for all $f\in \mathcal{F}_X$}\}= \mathcal{C}_X$
\end{itemize} Fr\"olicher and Kriegl \cite{fk}, and Kriegl and
Michor \cite{km} are our main reference for Fr\"olicher spaces.
The following terminology will be used in the paper: Given a
Fr\"olicher space $(X, \mathcal{C}_X, \mathcal{F}_X)$, the pair
$(\mathcal{C}_X, \mathcal{F}_X)$ is called a \textsl{smooth
structure}; the elements of $\mathcal{C}_X$ and $\mathcal{F}_X$
are called smooth curves and smooth functions respectively. The
topology assumed for a Fr\"olicher space $(X, \mathcal{C}_X,
\mathcal{F}_X)$ throughout the paper is the initial topology
$\mathcal{T}_\mathcal{F}$ induced by the set $\mathcal{F}_X$ of
functions. When there is no fear of confusion, a Fr\"olicher space
$(X, \mathcal{C}_X, \mathcal{F}_X)$ will simply be denoted $X$.
The most natural Fr\"olicher spaces are the finite dimensional
smooth manifolds, where if $X$ is such a smooth manifold, then
$\mathcal{C}_X$ and $\mathcal{F}_X$ consist of all smooth curves
$\mathbb{R}\rightarrow X$ and smooth functions $X\rightarrow
\mathbb{R}$. Euclidean finite dimensional smooth manifolds
$\mathbb{R}^n$, when viewed as Fr\"olicher spaces, are called
Euclidean Fr\"olicher spaces. In the sequel, by $\mathbb{R}^n$,
$n\in \mathbb{N}$, we mean the Fr\"olicher space $\mathbb{R}^n$,
equipped with its usual smooth manifold structure.

A Fr\"olicher space $X$ is called Hausdorff if and only if the
smooth real-valued functions on $X$ are point-separating, i.e. if
and only if $\mathcal{T}_\mathcal{F}$ is Hausdorff.

A Fr\"olicher structure $(\mathcal{C}_X, \mathcal{F}_X)$ on a set
$X$ is said to be generated by a set $F_0\subseteq \mathbb{R}^X$
(resp. $C_0\subseteq X^\mathbb{R}$) if $\mathcal{C}_X= \Gamma F_0$
and $\mathcal{F}_X= \Phi\Gamma F_0$ (resp. $\mathcal{F}_X= \Phi
C_0$ and $\mathcal{C}_X= \Gamma\Phi C_0$ ). Note that different
sets $F_0\subseteq \mathbb{R}^X$ on the same set $X$ may give rise
to a same smooth structure on $X$. A set mapping $\varphi:
X\rightarrow Y$ between Fr\"olicher spaces is called a map of
Fr\"olicher spaces or just a smooth map if for each $f\in
\mathcal{F}_Y$, the pull back $f\circ \varphi\in \mathcal{F}_X$.
This is equivalent to saying that for each $c\in \mathcal{C}_X$,
$\varphi\circ c\in \mathcal{C}_Y$. For Fr\"olicher spaces $X$ and
$Y$, $C^\infty(X, Y)$ will denote the collection of all the smooth
maps $X\rightarrow Y$. The resulting category of Fr\"olicher
spaces and smooth maps is denoted by $\mathbb{FRL}$.

Some useful facts regarding Fr\"olicher spaces can be gathered in
the following

\begin{theo} The category $\mathbb{FRL}$ is complete (i.e.
arbitrary limits exist ), cocomplete (i.e. arbitrary colimits
exist), and Cartesian closed. \end{theo}

Given a collection of Fr\"olicher spaces $\{X_i\}_{i\in I}$, let
$X= \prod_{i\in I}X_i$ be the set product of the sets
$\{X_i\}_{i\in I}$ and $\pi_i:X\rightarrow X_i$, $i\in I$, denote
the projection map $(x_i)_{i\in I}\mapsto x_i$. The initial
structure on $X$ is generated by the set \[F_0= \bigcup_{i\in
I}\{f\circ \pi_i:\ f\in \mathcal{F}_{X_i}\}.\]The ensuing
Fr\"olicher space $(X, \Gamma F_0, \varphi\Gamma F_0)$ is called
the product space of the family $\{X_i\}_{i\in I}$. Clearly,
\[\Gamma F_0= \{c: \mathbb{R}\rightarrow X|\ \mbox{if $c(t)=
(c_i(t))_{i\in I}$, then $c_i\in \mathcal{C}_{X_i}$ for every
$i\in I$}\}.\]Now, let $\biguplus_{i\in I}X_i$ be the disjoint
union of sets $\{X_i\}_{i\in I}$, and $\iota_{X_i}: X_i\rightarrow
\biguplus_{i\in I}X_i$ the inclusion map. Place the smooth final
structure on $\biguplus_{i\in I}X_i$ corresponding to the family
$\{\iota_{X_i}\}_{i\in I}$. The resulting Fr\"olicher space is
called the coproduct of $\{X_i\}_{i\in I}$, and denoted
$\coprod_{i\in I}X_i$, and \[\mathcal{F}_{\coprod_{i\in I}X_i}=
\{f: \coprod_{i\in I}X_i\rightarrow \mathbb{R}|\ \mbox{for each
$i\in I$, $f|_{X_i}\in \mathcal{F}_{X_i}$}\}\]is the collection of
smooth functions for the coproduct.

\begin{cor} Let $X$, $Y$, and $Z$ be Fr\"olicher spaces. Then the
following canonical mappings are smooth. \begin{itemize} \item ev:
$C^\infty(X, Y)\times X\rightarrow Y$, $ (f, x)\mapsto f(x)$ \item
ins:$X\rightarrow C^\infty(Y, X\times Y)$, $x\mapsto (y\mapsto
\mbox{ins$(x)(y)= (x, y)$})$ \item comp:$C^\infty(Y, Z)\times
C^\infty(X, Y)\rightarrow C^\infty(X, Z)$, $(g, f)\mapsto g\circ f$
\item $f_*: C^\infty(X, Y)\rightarrow C^\infty(X, Z)$, $f_*(g)=
f\circ g$, where $f\in C^\infty(Y, Z)$ \item $g_*: C^\infty(Z,
Y)\rightarrow C^\infty(X, Y)$, $g_*(f)= f\circ g$, where $g\in
C^\infty(X, Z)$.\end{itemize}\end{cor}

Given Fr\"olicher spaces $X$, $Y$, and $Z$; in view of the
cartesian closedness of the category $\mathbb{FRL}$, the
exponential law \[C^\infty(X\times Y, Z)\cong C^\infty(X,
C^\infty(Y, Z))\] holds. Because $\mathcal{F}_X= C^\infty(X,
\mathbb{R})$, it follows by cartesian closedness of $\mathbb{FRL}$
that the collection $\mathcal{F}_X$ can be made into a Fr\"olicher
space on its own right.

Finally we would like to show how to construct smooth braking
functions, following Hirsch \cite{hi}. Smooth braking functions
are tools that are behind most results in this paper. In
\cite{jn}, it is shown that the function
$\varphi:\mathbb{R}\rightarrow \mathbb{R}$ given by
\[\varphi (u)=\left\{\begin{array}{ll} 0 & \mbox{if $u\leq 0$}\\
e^{\frac{-1}{u}} & \mbox{if $u>0$}\end{array}\right.\]is smooth.
Substituting $x^2$ for $u$ in the above function, one sees that
the function $\psi:\mathbb{R}\rightarrow \mathbb{R}$, given by
\[\psi (x)=\left\{\begin{array}{ll} 0 & \mbox{if $x\leq 0$}\\
e^{\frac{-1}{x^2}} & \mbox{if $u>0$}\end{array}\right.\]is smooth.
Now, let us construct a smooth function $\alpha:
\mathbb{R}\rightarrow \mathbb{R}$ with the following properties.
Let $0\leq a<b$. $\alpha(t)$ should satisfy:
\begin{itemize} \item $\alpha(t)= 0$ for $t\leq a$, \\ \item
$0<\alpha(t)<1$ for $a<t<b$, \\ \item $\alpha$ is strictly
increasing for $a<t<b$,\\ \item $\alpha(t)= 1$ for $t\geq
b$.\end{itemize}Define $\alpha: \mathbb{R}\rightarrow [0, 1]$ by
\[\alpha(t)=
\frac{\int_a^t\gamma(x)dx}{\int_a^b\gamma(x)dx},\]where
$\gamma(x)= \psi(x-a)\psi(b-x)$.

In the sequel, the notation $\alpha_\epsilon$,
$0<\epsilon<\frac{1}{2}$, will refer to a smooth braking function
with the following properties \begin{itemize} \item
$\alpha_\epsilon(t)= 0$ for $t\leq \epsilon$, \\ \item $0<
\alpha_\epsilon(t)< 1$ for $\epsilon< t< 1-\epsilon$, \\ \item
$\alpha$ strictly increasing for $\epsilon< t< 1-\epsilon$, \\
\item $\alpha_\epsilon(t)= 1$ for $1-\epsilon \leq t$.
\end{itemize}

\section{Basic Constructions of Homotopy Theory in $\mathbb{FRL}$}

In this section, we define the fundamental notions of homotopy
theory in the category $\mathbb{FRL}$, such as the homotopy
relation and the mapping cylinder. We begin with an overview of
our approach to homotopy in $\mathbb{FRL}$, and then discuss
alternate Fr\"olicher structures on the unit interval which are
used in this and subsequent sections.

\subsection{Our Approach to Homotopy Theory in $\mathbb{FRL}$}
 One might begin investigating homotopy theory in
$\mathbb{FRL}$ by simply following the homotopy theory of
topological spaces, replacing continuous functions with smooth
ones. One can certainly define the notion of a homotopy $H:
I\times X\rightarrow Y$ between smooth maps $H(0, -)$ and $H(1,
-)$ in this way (which we do). One can even get as far as the left
Puppe sequence (see \cite{fch}), but eventually difficulties begin
to arise.

Extending functions defined on a subspace of a Fr\"olicher space
tends to be a little tricky, and so the definition of a
cofibration in $\mathbb{FRL}$ is one that needs careful
consideration. We envisage to construct the right Puppe sequence
in a future paper. To do this we define a slightly weaker notion
of cofibration than the notion obtained from topological spaces.
In addition, we define the mapping cylinder of a smooth map $f:
X\rightarrow Y$ using not the unit interval, but a modified
version called the {\bf weakly flattened unit interval}, denoted
$\mathbb{I}$, which, as one can show, is topologically
homeomorphic to the unit interval. This modified structure on the
unit interval allows us to show that the inclusion of a space $X$
into the mapping cylinder of $f: X\rightarrow Y$ is a cofibration
(in our weaker sense ).

The weakly flattened unit interval is useful, but it also has its
drawbacks. It would be ideal to have a single structure on the
unit interval that can be used throughout out homotopy theory, but
the weakly flattened unit interval is not suitable, because it has
the rather restrictive property that a smooth map $f: I\rightarrow
I$ on the usual unit interval often does not define a smooth map
$f: \mathbb{I}\rightarrow \mathbb{I}$ unless the endpoints of the
interval are mapped to the endpoints. This restrictive property
means that we only use the flattened unit intervals where they are
absolutely necessary.

In our future work, we will investigate whether with our modified
notions of cofibration and mapping cylinder, Baues' cofibration
axioms are satisfied.

\subsection{Flattened Structures on the Unit Interval}We define two main
Fr\"olicher structures which we call the {\bf flattened unit
interval } and the {weakly flattened unit interval }. Let
$(\mathcal{C}_I, \mathcal{F}_I)$ be the subspace structure induced
on $I$ by the inclusion $I\hookrightarrow \mathbb{R}$.

\begin{mydf} The Fr\"olicher space $(\mathbf{I},
\mathcal{C}_\mathbf{I}, \mathcal{F}_\mathbf{I})$, where the
structure $(\mathcal{C}_\mathbf{I}, \mathcal{F}_\mathbf{I})$ is
the structure generated by the set \begin{eqnarray*} F & = &
\{f\in \mathcal{F}_\mathbf{I}|\ \mbox{there exists $0<\epsilon
<\frac{1}{4}$ with $f(t)= f(0)$ for $t\in [0, \epsilon)$ and }\\
{} & {} & \mbox{$f(t)= f(1)$ for $t\in (1-\epsilon,
1]$}\},\end{eqnarray*}is called the flattened unit interval.
\end{mydf}

It is easy to see that any continuous map $c:
\mathbb{R}\rightarrow [0, 1]$ defines a structure curve on
$\mathbf{I}$ if and only if it is smooth at every point $t\in
\mathbb{R}$, where $c(t)\in (0, 1)$, .

We define the {\bf left} (resp. {\bf right}) {\bf flattened unit
interval}, denoted by ${\bf I}^-$ (resp. ${\bf I}^+$), to be the
Fr\"olicher space whose underlying set is the unit interval $[0,
1]$, and structure is the structure generated by the structure
functions in $\mathcal{F}_I$ that are constant near $0$ (resp.
$1$).

\begin{mydf} The Fr\"olicher space $(\mathbb{I},
\mathcal{C}_\mathbb{I}, \mathcal{F}_\mathbb{I})$, with the
structure defined below is called the {\bf weakly flattened unit
interval}. The underlying set is the unit interval; the structure
$(\mathcal{C}_\mathbb{I}, \mathcal{F}_\mathbb{I})$ is generated by
the family \[F= \{f\in \mathcal{F}_I|\ \lim_{t\rightarrow
0^+}\frac{d^n}{dt^n}f(t)= 0,\ \lim_{t\rightarrow
1^-}\frac{d^n}{dt^n}f(t)= 0,\ n\geq 1\}.\]We call the property,
for all $f\in F$,
\[\lim_{t\rightarrow 0^+}\frac{d^n}{dt^n}f(t)= 0,\
\lim_{t\rightarrow 1^-}\frac{d^n}{dt^n}f(t)= 0,\ n\geq 1,\]the
zero derivative property of $f$.\end{mydf}

We shall prove that all structure functions on $\mathbb{I}$ have
the zero derivative property, in other words,
$\mathcal{F}_\mathbb{I}= F$. To that effect, we need the following
lemma.

\begin{lem}\label{primary} Let $c: \mathbb{R}\rightarrow \mathbb{R}$ be a smooth
real-valued function at $t= t_0$, and let $f:
\mathbb{R}\rightarrow \mathbb{R}$ be a smooth real-valued function
at $t= c(t_0)$. Then,
\begin{eqnarray*} \frac{d^n}{dt^n}(f\circ
c)(t_0)  =  f^{(n)}(c(t_0))(c'(t_0))^n+ \mbox{terms of the form }\\
af^{(k)}(c(t_0))(c'(t_0))^{m_1}(c''(t_0))^{m_2}\ldots
(c^{(n-1)}(t_0))^{m_{n-1}},\end{eqnarray*} where $k<n$ and $a\in
\mathbb{R}$. In addition, if $a\neq 0$ then at least one of $m_2,
m_3, \ldots, m_{n-1}$ is also non-zero.\end{lem} \P The proof is
done by induction. For the sake of brevity, we call the term
$f^{(n)}(c(t_0))(c'(t_0))^n$ the primary term for $n$, and the
terms of the form
$af^{(k)}(c(t_0))(c'(t_0))^{m_1}(c''(t_0))^{m_2}\ldots
(c^{(n-1)}(t_0))^{m_{n-1}}$ the lower order terms for $n$. The
statement is true for $n=1$ and for $n=2$. Suppose the result is
true for $n=k$. To show that the result holds for $n=k+1$, since
\begin{eqnarray*} \frac{d^{k+1}}{dt^{k+1}}(f\circ c)(t_0)=
\frac{d}{dt}(f^{(k)}(c(t_0))(c'(t_0))^k) \\ + \mbox{terms of the
form
$\frac{d}{dt}(af^{(j)}(c(t_0))(c'(t_0))^{m_1}(c''(t_0))^{m_2}\ldots
(c^{(k-1)}(t_0))^{m_{k-1}}),$}
\end{eqnarray*}where $j<k+1$ and $a\in \mathbb{R}$, we need only show that
\[\frac{d}{dt}(af^{(j)}(c(t_0))(c'(t_0))^{m_1}(c''(t_0))^{m_2}\ldots
(c^{(k-1)}(t_0))^{m_{k-1}})\] gives rise to lower terms for $n=
k+1$, which is by the way straightforward. \proofterminator

\begin{theo} $\mathcal{F}_\mathbb{I}=\{f\in \mathcal{F}_I|\
\lim_{t\rightarrow 0^+}\frac{d^n}{dt^n}f(t)=0= \lim_{t\rightarrow
1^-}\frac{d^n}{dt^n}f(t)\}=:F$\end{theo} \P That $F\subseteq
\mathcal{F}_\mathbb{I}$ is evident. We must show the reverse
inequality. Let $0<\epsilon< \frac{1}{2}$, and $0<M<1$. Consider
the function $c_M:\mathbb{R}\rightarrow \mathbb{R}$, given by
\[c_M(t)= (1-\alpha_\epsilon(|t|))\beta_M(t)+
\alpha_\epsilon(|t|),\]where
$\alpha_\epsilon:\mathbb{R}\rightarrow \mathbb{R}$ is a smooth
braking function as defined in the Preliminaries, and
$\beta_M:\mathbb{R}\rightarrow \mathbb{R}$ is given by
\[\beta_M(t)= \left\{\begin{array}{ll} -Mt & \mbox{if $t\leq 0$}\\
t & \mbox{if $t>0$}\end{array}\right.\]It is easily seen that
$c_M$ is continuous over all $\mathbb{R}$, and smooth over all
$\mathbb{R}$ except at $t=0$. Also note that $0< c_M(t)<1$ for all
$t\in \mathbb{R}$, and $c_M(t)= \beta_M(t)=0$ for all $0\leq
t<\epsilon$. Now, \begin{eqnarray*}\frac{d}{dt}c_M(t)=
\frac{d}{dt}\beta_M(t) =-M, & & \mbox{for $-\epsilon< t< 0$} \\
\frac{d}{dt}c_M(t)= \frac{d}{dt}\beta_M(t)=1, & & \mbox{for
$0<t<\epsilon$}\end{eqnarray*}For $n>1$, we have
\[\frac{d^n}{dt^n}c_M(t)= \frac{d^n}{dt^n}\beta_M(t)= 0, \ \
\mbox{for $t\in (-\epsilon, 0)\cup (0, \epsilon)$}.\]We now show
that for $c_M\in \Gamma F$. To this end, let $f\in F$. To show
that $f\circ c_M:\mathbb{R}\rightarrow \mathbb{R}$ is smooth, it
is obvious that we need only concentrate on the point $t=0$,
because $f\circ c$ is smooth at every $t\neq 0$. It follows for
$t\neq 0$, and $n\in \mathbb{N}$ that Lemma \ref{primary} applies.
But as $t\rightarrow 0$, $c_M(t)\rightarrow 0^+$, and so, letting
$s=c_M(t)$, we have \[\lim_{t\rightarrow 0}f^{(j)}(c_M(t))=
\lim_{s\rightarrow 0^+}f^{(j)}(s)=0,\]for all $j\in \mathbb{N}$,
by the zero derivative property of $f$. Thus, as $t$ approaches
the value $0$, the primary term and all the lower order terms of
$\frac{d^n}{dt^n}(f\circ c_M)(t)$ vanish, and we have shown that
$f\circ c_M$ is smooth at $t=0$. This implies that $f\circ c_M\in
C^\infty(\mathbb{R}, \mathbb{R})$ for all $f\in F$. It follows
that $c_M\in \Gamma F$.

We are now ready to show that $\mathcal{F}_\mathbb{I}\subseteq F$.
To this end, suppose that we are given a structure function $f\in
\mathcal{F}_\mathbb{I}$. We shall show that this $f$ has the zero
derivative property, and is thus an element of $F$.

Since $f\in \mathcal{F}_\mathbb{I}$, we know that $f\circ c$ is a
smooth real-valued function for every $c\in \Gamma F$. In
particular, $f\circ c_M$ is smooth for all $0<M <1$. Thus, for any
$n\in \mathbb{N}$, \[\lim_{t\rightarrow
0^-}\frac{d^n}{dt^n}(f\circ c_M)(t)= \lim_{t\rightarrow
0^+}\frac{d^n}{dt^n}(f\circ c_M)(t).\]As $t\rightarrow 0^-$,
$c_M(t)\rightarrow 0^+$; let us consider the lower order terms for
$n$. Each term of the form
\[af^{(k)}(c_M(t))(c'_M(t))^{m_1}(c''_M(t))^{m_2}\ldots
(c_M^{(n-1)}(t))^{m_{n-1}}\]has some term $(c^{(i)}_M(t))^{m_i}$,
for some $i>1$, with $m_i\neq 0$. But $\lim_{t\rightarrow
0^-}c_M^{(i)}(t)=0$, if $i>1$, and so \[\lim_{t\rightarrow
0^-}af^{(k)}(c_M(t))(c'_M(t))^{m_1}(c''_M(t))^{m_2}\ldots
(c_M^{(n-1)}(t))^{m_{n-1}}=0.\]So all the lower order terms fall
away, therefore \begin{eqnarray*}\begin{array}{lll}
\lim_{t\rightarrow 0^-}\frac{d^n}{dt^n}(f\circ c_M)(t) & = &
\lim_{t\rightarrow 0^-}f^{(n)}(c_M(t))(c'_M(t))^n\\ {} & = &
\lim_{t\rightarrow 0^-}f^{(n)}(c_M(t))(-M)^n\\ {} & = &
\lim_{s\rightarrow
0^+}f^{(n)}(s)(-M)^n,\end{array}\end{eqnarray*}where $s= c_M(t)$.
In a similar way one shows that \[\lim_{t\rightarrow
0^+}\frac{d^n}{dt^n}(f\circ c_M)(t)= \lim_{s\rightarrow
0^+}f^{(n)}(s).\]But $f\circ c_M$ is smooth, therefore
$\lim_{s\rightarrow 0^+}f^{(n)}(s)(-M)^n= \lim_{s\rightarrow
0^+}f^{(n)}(s)$, which implies that $\lim_{s\rightarrow
0^+}f^{(n)}(s)=0.$

We have shown that the zero derivative property of $f$ holds for
the left endpoint of the unit interval. To show that the zero
derivative property of $f$ holds for the right endpoint of $f$,
note that $d_M: \mathbb{R}\rightarrow \mathbb{R}$, $d_M(t)=
1-c_M(t)$, is a smooth real-valued function with $d(0)=1$, and
$0\leq d_M(t) \leq 1$ for all $t\in \mathbb{R}$. One can follow a
similar procedure to the above, using $d_M$ instead of $c_M$ to
show that $\lim_{s\rightarrow 1^-}f^{(n)}=0.$\proofterminator

\subsection{Some Properties of Smooth Functions between the
Flattened Unit Intervals}One has to be careful when dealing with
the various flattened unit intervals. A smooth function $f:
I\rightarrow I$ from the $\mathbb{R}$- Fr\"olicher subspace unit
interval $I$ to itself need not define a smooth function $f:
\mathbf{I}\rightarrow \mathbf{I}$, for example. Conversely, not
every smooth function $f: \mathbf{I}\rightarrow \mathbf{I}$
defines a smooth function $f: I\rightarrow I$. In particular, we
need to be aware of the fact that addition and multiplication of
functions when defined between the various flattened unit
intervals does not preserve smoothness, as is the case with the
usual unit interval.

\begin{exam}\end{exam} The function $f: I\rightarrow I$, $f(t)=
\frac{1}{2}t$ is clearly smooth, but the corresponding function $f:
\mathbf{I}\rightarrow \mathbf{I}$, given by the same formula, is not
smooth. To see this, let $\alpha: \mathbb{R}\rightarrow \mathbb{R}$
be a smooth braking function with the properties that
\begin{itemize}\item $\alpha(t)= -1$, for $t<-\frac{3}{4}$, \item
$\alpha(t)=t$, for $-\frac{1}{4}< t< \frac{1}{4}$, \item
$\alpha(t)=1 $, for $t> \frac{3}{4}$.\end{itemize} Define $c:
\mathbb{R}\rightarrow \mathbf{I}$ by $c(t)= 1- |\alpha(t)|$. The
curve $c$  is smooth everywhere except at $t=0$, where $c(0)=1$.
However, every generating function $f$ on $\mathbf{I}$ is constant
near $1$, and so the composite $f\circ c$ is smooth. Thus $c$ is a
structure curve on $\mathbf{I}$. Now, $f\circ c:
\mathbb{R}\rightarrow \mathbf{I}$ is given by $(f\circ c)(t)=
\frac{1}{2}(1- |\alpha(t)|)$. Let $h: \mathbf{I}\rightarrow
\mathbb{R}$ be a structure function with the properties that
\begin{itemize} \item $h(s)= 0$, for $s<\frac{1}{8},$ \item $h(s)=
s$, for $\frac{1}{4}<s <\frac{3}{4}$, \item $h(s)= 1$, for
$\frac{7}{8}.$\end{itemize} Then $(h\circ f\circ c)(t)=
\frac{1}{2}(1- |\alpha(t)|)$ for $t$ near $0$, and is not smooth
at $t=0$. Thus $f$ does not define a smooth function from
$\mathbf{I}$ to $\mathbf{I}$.

\begin{exam}\end{exam} The function $f: \mathbf{I}\rightarrow \mathbf{I}$, $f(t)=
\sqrt{t}$, is smooth, but the corresponding $f: I\rightarrow I$,
given by the same formula, is not smooth. This follows from the fact
that $f$ is smooth on the open interval $(0, 1)$, and a generating
function $g$ on $\mathbf{I}$ is constant near $0$ and $1$. On the
side, $f:I\rightarrow I$ is not smooth, because if
$c:\mathbb{R}\rightarrow I$ is a structure curve with $c(t)= t^2$
near $t=0$, then $(f\circ c)(t)= |t|$ near $t=0$, which is not
smooth on $I$ at $t=0$.

\begin{exam}\end{exam} The functions $f, g: \mathbf{I}^-\rightarrow \mathbf{I}^-$,
given by $f(t)= \frac{1}{2}\sqrt{t}$ and $g(t)= \frac{1}{4}$ are
both smooth, but the sum $f(t)+ g(t)= \frac{1}{2}\sqrt{t}+
\frac{1}{4}$ is not smooth.

The following lemma follows from the definition of the Fr\"olicher
structures on the various flattened unit intervals.

\begin{lem} Let $f: I\rightarrow I$ be a smooth function with the
properties that $f(0)= 0$ and $f(1)= 1$. Then the following maps
are smooth: \begin{itemize} \item $f: I\rightarrow
\mathbf{I}^\pm,$ \item $f: I\rightarrow \mathbb{I},$ \item $f:
\mathbf{I}^\pm\rightarrow \mathbf{I},$ \item $f:
\mathbb{I}\rightarrow \mathbf{I},$ \item $f: I\rightarrow
\mathbf{I}.$\end{itemize}\end{lem}

The function defined in the following example is for later
reference.

\begin{exam}\label{445}\end{exam} Let $H: \mathbf{I}\times \mathbf{I}^-\rightarrow \mathbf{I}^-$
be given by $H(t, s)= (1-\alpha(t))s$, where $\alpha:
\mathbb{R}\rightarrow \mathbb{R}$ is a smooth braking function with
the properties that
\begin{itemize} \item $\alpha(t)=0$ for $t<\frac{1}{4},$ \item
$0\leq \alpha(t)\leq 1$ for all $t\in \mathbb{R},$ \item
$\alpha(t)=1$ for $t>\frac{3}{4}.$\end{itemize}We show that $H$ is
smooth. To see this, let $f: \mathbf{I}^-\rightarrow \mathbb{R}$
be a generating function on $\mathbf{I}^-$. So $f$ is constant
near $0$. Now, let $c: \mathbb{R}\rightarrow \mathbf{I}\times
\mathbf{I}^-$ be a structure curve, given by $c(v)= (t(v), s(v))$.
The curve $t$ is a structure curve on $\mathbf{I}$, and so is a
smooth real-valued function for all $v\in \mathbb{R}$, except
possibly when $t(v)=0$ or $t(v)=1$. Similarly, the curve $s$ is a
structure curve on $\mathbf{I}^-$, and so is smooth for all $v\in
\mathbb{R}$ except possibly when $s(v)=0$. Now consider the
composite $H\circ c: \mathbb{R}\rightarrow \mathbf{I}^-$. Clearly,
$\alpha(t(v))$ is smooth for all $v$, since the only possible
points for non-smoothness occur when $t(v)=0$ or $t(v)=1$, and
$\alpha(t(v))$ is locally constant near these points.
Consequently, $H\circ c$ is smooth everywhere except possibly when
$s(v)=0$. Now, let's consider $f\circ H\circ c:
\mathbb{R}\rightarrow \mathbb{R}$; the only possible points for
non-smoothness are those in which $s$ is $0$, i.e. $H\circ =0$.
But $f$ is a structure generating function on $\mathbf{I}^-$, and
so is locally constant near $0$. This shows that $f\circ H\circ c$
is smooth for all $v\in \mathbb{R}$, and thus $H$ is smooth.

\subsection{Homotopy in $\mathbb{FRL}$ and Related
Objects}\label{45}
\begin{mydf}
{\bf (1)} Let $X$ be a Fr\"olicher space, and $x_0$, $x_1\in X$.
We say that $x_0$ is smoothly path-connected to $x_1$ if there is
a smooth path $c:I\rightarrow X$ such that $c(0)= x_0$ and $c(1)=
x_1$. We write $x_0\simeq x_1$. The relation $\simeq$ is called
smooth homotopy when it is applied to hom-sets.

{\bf (2)} Let $f: X\rightarrow Y$ be a map of Fr\"olicher spaces.
$f$ is called a smooth homotopy equivalence provided there exists
a smooth map $g: Y\rightarrow X$ such that $f\circ g\simeq 1_Y$
and $g\circ f\simeq 1_X$.
\end{mydf}One can show that smooth homotopy is a congruence in
$\mathbb{RFL}$. In practice, we say that smooth maps $f, g:
X\rightarrow Y$ are smoothly homotopic if there exists a smooth map
$H: I\times X\rightarrow Y$ with $H(0, -)= f$ and $H(1, -)=g$. If
$A\subseteq X$ is subspace of $X$, then we say that $H$ is a smooth
homotopy (rel $A$) if the map $H$ has the additional property that
$H(t, a)=a$ for each $t\in I$ and $a\in A$. See Cherenack \cite{pch}
and Dugmore \cite{dug} for more detail regarding smooth homotopy.

The notion of deformation retract is fundamental to topological
homotopy theory. The following definitions are adapted for smooth
homotopy, and will be needed at a later stage.

\begin{mydf} Let $A\subseteq X$ be a subspace of a Fr\"olicher space
$X$, and let $i: A\hookrightarrow X$ denote the inclusion map. Then
\begin{itemize} \item We say that $A$ is a retract of $X$ if there
exists a smooth map $r:X\rightarrow A$ such that $ri= 1_A$. We call
$r$ a retraction. \item We call $A$ a weak deformation retract of
$X$ if the inclusion $i$ is a smooth homotopy equivalence. \item The
subspace $A$ is called a deformation retract of $X$ if there exists
a retraction $r: X\rightarrow A$ such that $ir\simeq 1_X$. \item The
subspace $A$ is called a strong deformation retract of $X$ if there
exists a retraction $r: X\rightarrow A$ such that $ir\simeq 1_X(rel
A)$.\end{itemize}\end{mydf}

\begin{mydf} The mapping cylinder $I_f$ of $f:X\rightarrow Y$ is
defined by the following pushout \[\xymatrix{X\ar[r]^f\ar[d]^{i_1} & Y\ar[d] \\
I\times X\ar[r] & I_f}\]where $i_1: X\rightarrow I\times X$ is given
by $i_1(x)= (1, x)$, for any $x\in X$. We denote the elements of
$I_f$ by $[t, x]$ or $[y]$, where $(t, x)\in I\times X$ and $y\in
Y$.

Replacing $I\times X$ in the above pushout diagram by
$\mathbf{I}\times X$ or $\mathbb{I}\times X$, we obtain the
flattened mapping cylinder $\mathbf{I}_f$ and weakly flattened
mapping cylinder $\mathbb{I}_f$ of $f$ respectively. We use the same
notation for elements of these flattened mapping cylinders as
described above for the mapping cylinder.
\end{mydf}There is also a map $i_0: X\rightarrow I\times X$, defined
by $i_0(x)= (0, x)$ for $x\in X$. This induces an inclusion map
$i'_0: X\rightarrow I_f$, which identifies $X$ with the
Fr\"olicher subspace $i'_0(X)$ of $I_f$. An inclusion is induced
in a similar way for the flattened mapping cylinders. If one
identifies $\{0\}\times X$ to a point in the mapping cylinder
$I_f$ of a map $f: X\rightarrow Y$, then one obtains the
$\mathbf{mapping\  cone}$ $T_f$ of the map $f$. In a similar
fashion, we define the $\mathbf{flattened\ mapping\ cone}$
$\mathbf{T}_f$ and $\mathbf{weakly\ flattened\ mapping\ cone}$
$\mathbb{T}_f$ of a smooth map $f: X\rightarrow Y$.

\subsection{Cofibrations in $\mathbb{FRL}$} A cofibration is a map
$i: A\rightarrow X$ for which the problem of extending functions
from $i(A)$ to $X$ is a homotopy problem. In other words, if a map
$f: i(A)\rightarrow Z$ can be extended to a map $f^*: X\rightarrow
Z$, then so can any map homotopic to $f$. For topological spaces,
the usual definition is phrased in a slightly more restrictive way.
The extension of a map $g\simeq_Hf$, for some homotopy $H: I\times
i(A)\rightarrow Z$, is required to exist at every level of the
homotopy simultaneously. In other words, one requires each $H(t, -)$
to be extendable in such a way that the resulting homotopy $H^*:
I\times X\rightarrow Z$ is continuous.

We weaken this definition somewhat, to enable smooth homotopy
extensions to be more easily constructed using a flattening at the
endpoints of the homotopy. This enables us to characterize smooth
cofibrations in terms of a flattened unit interval, and then later
to relate smooth cofibrations to smooth neighborhood deformation
retracts. Our definition of smooth cofibration, though different
from from Cap's definition, see \cite{cap}, leads to several
classical results as does Cap's. As pointed out by Cap, the analogue
of the classical definition of cofibration would not allow even
$\{0\}\hookrightarrow I$ to be a smooth cofibration. So, we have the
following

\begin{mydf} \label{461}A smooth map $i:A\rightarrow X$ is called a
smooth cofibration if, corresponding to to every commutative diagram
of the form \[\xymatrix{A \ar[r]^i\ar[d]_{(0, 1_A)} & X\ar[r]^f & Z\\
I\times A\ar[rru]_G},\]there exists a commutative diagram in
$\mathbb{FRL}$ of the form \[\xymatrix{X\ar[r]^f\ar[d]_{(0, 1_X)} &
Z\\  I\times X \ar[ur]^F & I \times A \ar[u]_{G'}\ar[l]^{1_I\times
i}},\]where $G': I\times A\rightarrow Z$ is given by $G'(t, a)=
G(\alpha_\epsilon(t), a)$ for some $0< \epsilon <\frac{1}{2}$, and
each $t\in I$, $a\in A$.
\end{mydf}

The problem of extending a map smoothly from a subspace of a
Fr\"olicher space to the whole space is a more difficult problem
than simply extending continuously. It is mainly for this reason
that the definition of smooth cofibration differs somewhat from
the corresponding definition of a topological cofibration.

\begin{lem} Let $i: A\rightarrow X$ be a smooth cofibration, then
$i$ is an initial morphism in $\mathbb{FRL}$. In addition, if $A$ is
Hausdorff, then $i$ is injective. So in this case $A$ can be
regarded as a subspace of $X$.\end{lem} \P Let us show that every
smooth map $f: A\rightarrow \mathbb{R}$ factors through $i$, that is
for every $f\in \mathcal{F}_A$, there exists $\tilde{f}\in
\mathcal{F}_X$ such that $f= \tilde{f}\circ i$. To this end,
consider the smooth map $G: I\times A\rightarrow \mathbb{R}$, given
by $H(t, a)= tf(a)$. Clearly, $0|_A= G(0, -)$, where $0:X\rightarrow
\mathbb{R}$ is the constant map $0$. It follows that there is map
$F: I\times X\rightarrow \mathbb{R}$ such that $F\circ (1\times i)=
G'$. Then, clearly $\tilde{f}:= F(1, -)$ has the desired property.

The remaining part of the proof of Proposition 3.3, in \cite{cap},
holds verbatim here as well. \proofterminator

In this paper, we are interested only in cofibrations that are
injective. Henceforth, all cofibrations are assumed to be injective.

All topological cofibrations are inclusions, and this result is
true for smooth cofibrations too. The proof of the following lemma
is essentially the same as the proof given by James \cite{ja} for
the topological result, although James's proof is in some sense
dual to ours, using path-spaces in place of cartesian products and
the adjoint versions of our homotopies.

\begin{lem} A cofibration \[\xymatrix{A\hspace{2mm}\ar@{>->}[r]^i & X}\]is a
smooth inclusion. \end{lem} \P Let $\mathbb{I}_i$ be a mapping
cylinder of $i$, and let $j: X\rightarrow \mathbb{I}_i$ be the
standard inclusion map. Consider the smooth map
$\gamma:I\rightarrow I$, $\gamma(t)= 1-t$, for all $t\in I$, and
the quotient map $q: (I\times A)\sqcup X\rightarrow \mathbb{I}_i$;
we have the following commutative diagram
\[\xymatrix{A \hspace{2mm}\ar@{>->}[r]^i\ar[d]_{(0, 1_A)} & X\ar[r]^j & \mathbb{I}_i\\
I\times A\ar[d]_{\gamma\times 1_A}\ar[rru]_q \\ I\times
A\ar[uurr]_G},\]where $G(t, a)= [(1-t, a)]$. Notice that the map
$G$ is smooth. Since $i$ is a cofibration, we have the
commutative diagram \[\xymatrix{X\ar[r]^j\ar[d]_{(0, 1_X)} & \mathbb{I}_i\\
I\times X \ar[ur]^F & I \times A \ar[u]_{G'}\ar[l]^{1_I\times
i}},\]where $G'(t, a)= G(\alpha_\epsilon(t), a)$ for some
$0<\epsilon <\frac{1}{2}$. Define $U: X\rightarrow \mathbb{I}_i$
by $U(x)= F(1, x)$. We have $U\circ i= G'(1, -)$, where $G'(1, a)=
[(0, a)]$, for every $a\in A$. Thus the assignment $a\mapsto G'(1,
a)$ defines the usual inclusion of $A$ into the mapping cylinder.
From this we deduce that $U\circ i$ is an inclusion, and hence $i$
is an inclusion. \proofterminator

There is an equivalent formulation of definition \ref{461}, given
in the following lemma.

\begin{lem} \label{463}A smooth map \[\xymatrix{A\hspace{2mm}\ar@{>->}[r]^i & X}\]
is a cofibration if and only if, for every smooth map $h: (0\times
X)\cup (\mathbf{I}^-\times i(A))\rightarrow Z$, the following
diagram \[\xymatrix{(0\times X)\cup (\mathbf{I}^-\times
i(A))\ar[r]^{\hspace{15mm}h}\ar[d]_j & Z\\ \mathbf{I}^-\times
X\ar[ru]_G},\]where $j$ is the evident inclusion, exists in
$\mathbb{FRL}$.
\end{lem} \P Suppose that the inclusion $\xymatrix{A\hspace{2mm}\ar@{>->}[r]^i &
X}$ is a smooth cofibration, and suppose that $h:(0\times X)\cup
(\mathbf{I}^-\times i(A))\rightarrow Z$ is a smooth map. We have the
diagram \[\xymatrix{(0\times B)\cup (\mathbf{I}^-\times
i(A))\ar[r]^{\hspace{15mm}h}\ar[d]_j & Z\\ \mathbf{I}^-\times
X}.\]We need to fill in a smooth map $G:\mathbf{I}^-\times
X\rightarrow Z$ which makes the resulting diagram commute. To do
this, notice that $h|_{\mathbf{I}^-\times i(A)}$ is smooth, and thus
the corresponding map $h|{I\times i(A)}$, using the usual unit
interval, is also smooth. We have the following diagram
\[\xymatrix{A\hspace{2mm}\ar@{>->}[r]^i\ar[d]_{(0, 1_A)} & X
\ar[r]^{h|_{0\times X}} & Z\\ I\times
A\ar[rru]_{h_{\mathbf{I}^-\times A}}},\]where $h|_{0\times X}(0, -):
X\rightarrow Z$ is denoted as $h|_{0\times X}$. The fact that $i$ is
a smooth cofibration yields the following $\mathbb{FRL}$-commutative
diagram: \[\xymatrix{X\ar[r]^{h_{0\times X}}\ar[d]_{(0, 1_A)} & Z\\
I\times X\ar[ru]^F & I\times A \ar[l]^{1_I\times
i}\ar[u]_{(h|_{\mathbf{I}^-\times A})'}},\]where
$(h|_{\mathbf{I}^-\times A})'(t, a)= h|_{\mathbf{I}^-\times
A}(\alpha_\epsilon(t), a)$, for some $0<\epsilon <\frac{1}{2}$. Now,
chose a smooth braking function $\beta: \mathbb{R}\rightarrow
\mathbb{R}$ with the following properties. \begin{itemize}\item
$\alpha(t)= 0$ for $t<\frac{\epsilon}{2}$, \item $\alpha(t)= t$ for
$\epsilon< t$.\end{itemize} $F$ may not be smooth on
$\mathbf{I}^-\times A$ due to the flattening requirements of the
left flattened unit interval. To correct this, set $G(t, a)=
F(\beta(t), a)$. Notice that the insertion of this braking function
does not affect the commutativity conditions of $G$, since the only
adjustments to $F$ occur in the first coordinate where the map
$(h|_{\mathbf{I}^-\times X})'$ is constant.

Now, assume the converse, i.e. to every smooth map $h: (0\times
X)\cup (\mathbf{I}^-\times i(A))\rightarrow Z$, corresponds a
commutative diagram \[\xymatrix{(0\times X)\cup (\mathbf{I}^-\times
i(A))\ar[r]^{\hspace{15mm}h}\ar[d]_j & Z\\ \mathbf{I}^-\times
X\ar[ru]_H}.\] We wish to show that the inclusion $i: A\rightarrow
X$ is a cofibration; so assume we have the following
diagram \[\xymatrix{A\ar[r]^i\ar[d]_{(0, 1_A)} & X \ar[r]^f & Z\\
I\times A \ar[rru]_G}.\]There exists the diagram \[\xymatrix{A\ar[r]^i\ar[d]_{(0, 1_A)} & X \ar[r]^f & Z\\
I\times A \ar[rru]_{G'}},\]where $G'(t, a)= G(\alpha_\epsilon(t),
a)$. Our hypothesis allows us to construct the diagram
\[\xymatrix{(0\times X)\cup (\mathbf{I}^-\times
i(A))\ar[r]^{\hspace{15mm}f\cup G'}\ar[d]_j & Z\\ \mathbf{I}^-\times
X\ar[ru]_H}.\]Note that $f\cup G'$ is smooth since
$\alpha_\epsilon(t)$ is constant near $0$. Since $H$ is smooth on
$\mathbf{I}^-\times X$ it defines a smooth map on $I\times X$. One
can verify that the diagram \[\xymatrix{X\ar[r]^f\ar[d]_{(0, 1_X)} &
Z\\ I\times X\ar[ru]^H & I\times A\ar[l]^{1_I\times
i}\ar[u]_{G'}}\]commutes as required. \proofterminator

\section{Smooth Neighborhood Deformation Retracts} This section is
concerned with the formulation of a suitable notion of smooth
neighborhood deformation retract. For topological spaces, the
statement that a closed subspace $A$ of $X$ is a neighborhood
deformation retract of $X$ is equivalent to the statement that the
inclusion $i: A\hookrightarrow X$ is a closed cofibration. We show
that in the category of Fr\"olicher spaces there is a notion of
smooth neighborhood deformation retract that gives rise to an
analogous result that a closed Fr\"olicher subspace $A$ of the
Fr\"olicher space $X$ is a smooth neighborhood deformation retract
of $X$ if and only if the inclusion $i: A\hookrightarrow X$ comes
from a certain subclass of cofibrations. As an application, we
construct the right Puppe sequence.

\subsection{SNDR pairs and SDR pairs} The definition of `smooth
neighborhood deformation retract' that we adopt in this paper is
similar to the definition of `R-SNDR pair'suggested in \cite{dug},
but we have modified the definition in order to retain only the
essential aspects of `first coordinate independence' defined in
\cite{dug}.

We begin by defining the `first coordinate independence property' of
a function on a product of a Fr\"olicher space with $\mathbf{I}$ (or
$\mathbf{I}^-$, $\mathbf{I}^+$).

\begin{mydf} Let $i: A\rightarrow X$ be a smooth map, and $c:
\mathbb{R}\rightarrow X$ a structure curve on $X$. Define
\begin{eqnarray*}\Lambda(c, i) & = &  \{t_*\in c^{-1}(i(A))|\ \mbox{there exists a sequence
$\{t_n\}$ of real numbers} \\ & & \mbox{with $\lim_{n\rightarrow
\infty}t_n= t_*$ and each $t_n\in c^{-1}(X-i(A))$}\}.\end{eqnarray*}
\end{mydf} The points in $\Lambda(c, i)$ are those values in
$\mathbb{R}$ where the curve `enters' $i(A)$ from $X-i(A)$, or
`touches' a point in $i(A)$ whilst remaining in $X-i(A)$ nearby.
Now, we are ready to define the `first coordinate independence
property' for a structure function on a product.

\begin{mydf}\label{512} Let $i: A\rightarrow X$ be a smooth map and suppose $f:
\mathbf{I}\times X\rightarrow \mathbb{R}$ is a structure function on
$\mathbf{I}\times X$. Let $c:\mathbb{R}\rightarrow \mathbf{I}\times
X$, given by $c(s)= (t(s), x(s))$ have the following properties
\begin{itemize}\item The map $x(s)$ is a structure curve on $X$.
\item For all $\epsilon >0$, $t(s)$ is a smooth real-valued
function on $\mathbb{R}-\cup_{s_*\in \Lambda(x, i)}[s_*-\epsilon,
s_*+\epsilon]$.\end{itemize}If, for every such map $c$, the
composite $f\circ c$ is a smooth real-valued function, then we say
that $f: \mathbf{I}\times X\rightarrow \mathbb{R}$ has the first
independence property (FCIP) with respect to $i$.

Extending the definition, we say that a map $g: \mathbf{I}\times
X\rightarrow Y$ has the FCIP with respect to $i$ if the composite
$h\circ g: \mathbf{I}\times X\rightarrow \mathbb{R}$ has the FCIP
with respect to $i$ for every $h\in \mathcal{F}_Y$.
\end{mydf}

Notice that we can formulate a similar definition of the FCIP if
we replace $\mathbf{I}$ throughout by $\mathbf{I}^-$ or
$\mathbf{I}^+$, leaving the rest of the definition unchanged. We
will have occasion to use this type of first coordinate
independence property in the later part of this work.

\vspace{5mm}\noindent\textbf{Note.} Let $i: A\rightarrow X$, and
suppose that we are given a map $g: \mathbf{I}\times X\rightarrow
Y$. Let $f: Y\rightarrow \mathbb{R}$ be a structure function on $Y$,
and suppose that $f\circ g: \mathbf{I}\times X\rightarrow
\mathbb{R}$ has the FCIP with respect to $i$ for any such $f$. Then,
given a smooth map $h: Y\rightarrow Z$, the composite $f'\circ
h\circ g: \mathbf{I}\times X\rightarrow \mathbb{R}$ has the FCIP
with respect to $i$ for any structure function $f'$ on $Z$.

The above note applies equally well if $g: \mathbf{I}^-\times
X\rightarrow Y$ or $g: \mathbf{I}^+\times X\rightarrow Y$ has the
FCIP with respect to $i$ when composed with a smooth function $h$ on
$Y$.

\begin{exam}\label{514}\end{exam}1. For any $i: A\rightarrow X$, the
projection onto the second coordinate $\pi_X: \mathbf{I}\times
X\rightarrow X$ has the FCIP.

\noindent 2. Let $\alpha: \mathbb{R}\rightarrow \mathbb{R}$ be a
smooth braking function with the properties that \begin{itemize}
\item $\alpha(t)=0$ if $t<\frac{1}{4}$, \item $0< \alpha(t)< 1$ if
$\frac{1}{4}\leq t\leq \frac{3}{4}$, \item $\alpha(t)= 1$ if
$\frac{3}{4}< t$. \end{itemize} Consider $0\hookrightarrow
\mathbf{I}^-$. Let $H: \mathbf{I}\times \mathbf{I}^-\rightarrow
\mathbf{I}^-$ be given by $H(t, s)= (1-\alpha(t))s$. Then, $f\circ
H: \mathbf{I}\times \mathbf{I}^-\rightarrow \mathbb{R}$ has the FCIP
with respect to the inclusion $0\hookrightarrow \mathbf{I}^-$, for
any $f\in \mathcal{F}_{\mathbf{I}^-}$.

\begin{mydf} Consider a smooth inclusion $i: A\hookrightarrow X$.
Suppose that there exists a smooth map $u: X\rightarrow
\mathbf{I}$, with $u^{-1}(0)= i(A)$. If there exists a smooth map
$H: \mathbf{I}\times X\rightarrow X$ that satisfies the following
properties: \begin{itemize} \item $H$ has the FCIP with respect to
$i$. \item $H(0, x)= x$ for all $x\in X$. \item $H(t, x)= x$ for
all $(t, x)\in \mathbf{I}\times i(A)$. \item $H(1, x)\in i(A)$ for
all $x\in X$ with $u(x)<1$,\end{itemize} then the pair $(X, A)$ is
called a smooth neighborhood deformation retract pair, or SNDR
pair for short.

If, in addition, $H$ is such that $H(1\times X)\subset i(A)$, then
the pair $(X, A)$ is called a smooth deformation retract pair, or an
SDR pair for short.

The subspace $A$ is called a smooth neighborhood deformation retract
or smooth deformation retract of $X$ if $(X, A)$ is an SNDR pair or
SDR pair, respectively.

The pair $(u, H)$ is called a representation for the SNDR (or SDR)
pair. \end{mydf}

\begin{exam}\end{exam}1. The pair $(X, \emptyset)$ is an SNDR pair. A representation
is $u(x)=1$, $H(t, x)= x$, for each $t\in \mathbf{I}$ and $x\in X$.

\noindent 2. The pair $(X, X)$ is an SNDR pair. A representation is
$u(X)=0$, $H(t, x)=x$, for each $t\in \mathbf{I}$ and $x\in X$.

\begin{lem}\label{517} The pair $(\mathbf{I}^-, 0)$ is an SDR pair.\end{lem} \P
Let $\alpha: \mathbb{R}\rightarrow \mathbb{R}$ be the smooth braking
function of Examples \ref{514}. A representation for $(\mathbf{I}^-,
0)$ as an SDR pair is $(u, H)$, where $u: \mathbf{I}^-\rightarrow
\mathbf{I}$ and $H: \mathbf{I}\times \mathbf{I}^-\rightarrow
\mathbf{I}^-$ are given by $u(s)=s$, and $H(t, s)= (1-\alpha(t))s$.
Clearly, the identity $u: \mathbf{I}^-\rightarrow \mathbf{I}$ is
smooth. And the map $H$, as shown in Example \ref{445}, is smooth
and clearly has the FCIP with respect to the inclusion, since
whenever $v$ approaches a value for which $s(v)= 0$, one has
\[g((1-\alpha(t(v)))s(v))= g(0)\]for $v$ in a neighborhood of this
value and $g\in \mathcal{F}_{\mathbf{I}^-}$. \proofterminator

\begin{lem} The pair $(\mathbf{I}, \{0, 1\})$ is an SNDR pair.
\end{lem} \P A representation $(u, H)$ for the SNDR pair can be
given as follows. Define $u: \mathbf{I}\rightarrow \mathbf{I}$ to be
a bump function such that \begin{itemize}\item $u(t)=0$ for $t=0$ or
$t=1$, \item $u(t)=1$ for $t\in [\frac{1}{4}, \frac{3}{4}]$, \item
$0<u(t) <1$ otherwise,\end{itemize}and let $\beta:
\mathbf{I}\rightarrow \mathbf{I}$ be a braking function with the
properties that $\beta(s)=0$ for $0\leq s \leq \frac{1}{4}$, and
$\beta(s)=1$ for $\frac{3}{4}\leq s\leq 1$. Let $0<\epsilon
\frac{1}{4}$, and define $H: \mathbf{I}\times \mathbf{I}\rightarrow
\mathbf{I}$ by $H(t, s)= (1-\alpha_\epsilon(t))s+
\alpha_\epsilon(t)\beta(s)$. It is clear that $H(0, s)= s$, $H(t,
0)= 0$, and $H(t, 1)= 1$. Suppose that $u(s)< 1$. Then, $s\in [0,
\frac{1}{4})\cup (\frac{3}{4}, 1]$. This implies that $\beta(s)=0$
or $\beta(s)= 1$. We then have $H(1, s)=0$ or $H(1, s)= 1$, which
means that $H(1, s)\in \{0, 1\}$ if $u(s)<1$.

To see that $H$ is smooth, let $f: \mathbf{I}\rightarrow \mathbb{R}$
be a generating function for the flattened unit interval. The only
possible points of non-smoothness are points where $t= 0, 1$ and $s=
0, 1$. The braking function $\alpha_\epsilon$ ensures that $H$ is
locally constant in the $t$b variable whenever $t$ is near $0$ or
$1$, so no problem arises from the $t$ component. When $s$ is near
$s=0$, we have $H(t, s)$ near $0$, and so the generating function
$f$ is locally constant. Similarly, when $s$ is near $s=1$, we have
$H(t, s)$ near $1$, and the generating function $f$ is again locally
constant. \proofterminator

We now show that the product of SNDR pairs is again an SNDR pair.

\begin{theo}\label{519}Let $i: A\hookrightarrow X$ and $j: B\hookrightarrow Y$
be inclusion mappings. If $(X, A)$ and $(Y, B)$ are SNDR pairs,
then so is \[(X\times Y, (X\times B)\cup (A\times Y)).\]If one of
$(X, A)$ or $(Y, B)$ is an SDR pair, then so is the pair
\[(X\times Y, (X\times B)\cup (A\times Y)).\]\end{theo}\P Let
$\alpha: \mathbb{R}\rightarrow I$ be a smooth braking function
with the properties that $\alpha(t)=0$ for $t\leq \frac{1}{4}$,
and $\alpha(t)=1$ for $t\geq \frac{3}{4}$, and let
$\beta:\mathbb{R}\rightarrow \mathbb{R}$ be a smooth increasing
braking function with the properties that $\beta(t)= t$ for $t\leq
\frac{1}{4}$, and $\beta(t)= 1$ for $t\geq \frac{3}{4}$. Suppose
that $(u, H)$ and $(v, J)$ are representations for the SNDR pairs
$(X, A)$ and $(Y, B)$, respectively. Let $\overline{u}:
X\rightarrow \mathbf{I}$, and $\overline{v}: Y\rightarrow
\mathbf{I}$ be given by $\overline{u}(x)= \beta(u(x))$ and
$\overline{v}(y)= \beta(v(y))$ respectively. Define $w: X\times
Y\rightarrow \mathbf{I}$ by $w(x, y)=
\overline{u}(x)\overline{v}(y)$. The braking function $\beta$
ensures smoothness of $\overline{u}$ and $\overline{v}$, and
consequently of $w$. We have $w^{-1}(0)= (X\times B)\cup (A\times
Y)$, as required. Define $Q: \mathbf{I}\times X\times Y\rightarrow
X\times Y$ as follows . \[Q(t, x, y)= \left\{\begin{array}{ll}
(H(\alpha(t), x), J(\alpha(t), y)) & \mbox{if $u(x)= v(y)=0$}\\
(H(\alpha(t), x),
J(\alpha(\frac{\overline{u}(x)}{\overline{v}(y)})\alpha(t), y)) &
\mbox{if $\overline{v}(y)\geq \overline{u}(x)$, $\overline{v}(y)>
0$,}\\ (H(\alpha(\frac{\overline{v}(y)}{\overline{u}(x)})\alpha(t),
x), J(\alpha(t), y)) & \mbox{if $\overline{u}(x)\geq
\overline{v}(y)$, $\overline{u}(x)> 0$.}\end{array}\right.\]We must
show that $Q$ is a smooth map, with the first coordinate
independence property with respect to the inclusion $(X\times B)\cup
(A\times Y)\hookrightarrow X\times Y$. We first consider each part
of the definition of $Q$ separately. The first part is clearly
smooth. Let us verify that $Q$ is smooth on the second part of its
definition; the third part is similar.

We need only focus on the component
$J(\alpha(\frac{\overline{u}(x)}{\overline{v}(y)})\alpha(t), y)$.
Each function making up
$J(\alpha(\frac{\overline{u}(x)}{\overline{v}(y)})\alpha(t), y)$ is
smooth individually, so we need only pay extra attention to those
parts that involve flattened unit intervals, remembering that
addition and multiplication on the flattened unit interval need not
preserve smoothness, as is the case for the usual unit interval.

So let us consider
$\alpha(\frac{\overline{u}(x)}{\overline{v}(y)})$; it is smooth
except possibly when $\frac{\overline{u}(x)}{\overline{v}(y)}$
approaches $0$ or $1$, since it is here that structure curves on the
flattened unit interval need not be smooth in the usual sense.
Clearly, if $u(x)$ approaches $0$ and $v(y)$ does not approach $0$,
then the braking function $\alpha$ ensures that
$\frac{\overline{u}(x)}{\overline{v}(y)}=0$ near such points. If
$\overline{v}(y)$ approaches $0$, then $\overline{u}(x)$ must
approach $0$ too. This situation is dealt with later.

Thus, $Q$, in part two of the definition, is smooth, and one can
show similarly that $Q$ in the third part of the definition is
smooth as well.

Let us now consider the overlaps of the three parts of the
definition of $Q$. Observe that if $\overline{u}(x)$ is in a
sufficiently small neighborhood of $\overline{v}(y)$, with
$\overline{u}(x)\neq 0$ and $\overline{v}(y)\neq 0$, then we have
$\alpha(\frac{\overline{u}(x)}{\overline{v}(y)}=
\alpha(\frac{\overline{v}(y)}{\overline{u}(x)})=1$, and so the
second and third parts of the definition of $Q$ coincide here. Thus,
it remains only to show that $Q$ is smooth as $\overline{u}(x)$ and
$\overline{v}(y)$ both approach $0$.

If $Q$ is smooth in each of its coordinates then it is smooth, so
consider the coordinate involving the map $J$. Let $c:
\mathbb{R}\rightarrow \mathbf{I}\times X\times Y$ be a structure
that is given by $c(s)= (t(s), x(s), y(s))$. Then, the map $c_1:
\mathbb{R}\rightarrow \mathbf{I}\times Y$, given by \[c_1(s)=
\left\{\begin{array}{ll} (\alpha(t(s)), y(s)) & \mbox{if
$\overline{u}(x(s))= \overline{v}(y(s))=0$}\\
(\alpha(\frac{\overline{u}(x(s))}{\overline{v}(y(s))})\alpha(t(s)),
y(s)) & \mbox{if $\overline{v}(y(s))\geq \overline{u}(x(s))$,
$\overline{v}(y(s))>0$}\\ (\alpha(t(s)), y(s)) & \mbox{if
$\overline{u}(x(s))\geq \overline{v}(y(s))$,
$\overline{u}(x(s))>0$}\end{array}\right.\]is a map satisfying the
conditions of Definition \ref{512}, since its second coordinate is
smooth, but its first coordinate may be singular as
$\overline{v}(y(s))$ ( and hence $\overline{u}(x(s))$) approaches
$0$. Since $J$ has the first coordinate independence property, the
map \[(Joc_1)(s)= \left\{\begin{array}{ll} J(\alpha(t(s)), y(s)) &
\mbox{if
$\overline{u}(x(s))= \overline{v}(y(s))=0$}\\
J(\alpha(\frac{\overline{u}(x(s))}{\overline{v}(y(s))})\alpha(t(s)),
y(s)) & \mbox{if $\overline{v}(y(s))\geq \overline{u}(x(s))$,
$\overline{v}(y(s))>0$}\\ J(\alpha(t(s)), y(s)) & \mbox{if
$\overline{u}(x(s))\geq \overline{v}(y(s))$,
$\overline{u}(x(s))>0$}\end{array}\right.\]is smooth. Thus,
$Q\circ c$ is smooth, and since $c$ is arbitrary, $Q$ is smooth.
In a similar way, the coordinate of $Q$ involving $H$ can be shown
to be smooth.

We now verify that $Q$ satisfies the required boundary conditions.
When $t=0$, all three lines defining $Q$ reduce to $(H(0, x), J(0,
y))=(x, y)$. Let $x\in A$ and $y\in B$; then $\overline{u}(x)=
\overline{v}(y)=0$. Therefore, $Q$ reduces to $(H(\alpha(t), x),
J(\alpha(t), y))= (x, y)$. If $x\in A$ and $y\notin B$, then $Q$
is given by the second part of its definition, which reduces to
$(H(\alpha(t), x), J(0, y))$. The case when $x\notin A$ and $y\in
B$ is similar. If $t=1$ and $0< w(x, y)< 1$ then either $0<
\overline{u}(x)< 1$ or $0< \overline{v}(y)< 1$. Suppose that $0<
\overline{u}(x)< 1$. Then either $\overline{u}(x)\leq
\overline{v}(y)$ or $\overline{v}(y)< \overline{u}(x)$. If
$\overline{u}(x)\leq \overline{v}(y)$, then $Q$ is given by the
second part of its definition, which reduces to $(H(1, x),
J(\alpha(\frac{\overline{u}(x)}{\overline{v}(y)}, y))\in
i(A)\times Y$. If $\overline{v}(y)< \overline{u}(x)$, then the
third part of the definition of $Q$ applies and $Q$ reduces to
$(H(\alpha(\frac{\overline{v}(y)}{\overline{u}(x)}), x), J(1,
y))\in X\times j(B)$.

Finally, we must show that for any $f\in \mathcal{F}_{X\times Y}$,
$f\circ Q$ has the first coordinate independence property with
respect to the inclusion $(X\times B)\cup (A\times
Y)\hookrightarrow X\times Y$. To this end, consider a map $c:
\mathbb{R}\rightarrow \mathbf{I}\times X\times Y$, given by $c(s)=
(t(s), x(s), y(s))$. Let $\{s_n\}$ be a sequence of real numbers
converging to $s_*$ with $c(s_n)\in (X\times Y)-((A\times Y)\cup
(X\times B))$, and $c(s_*)\in (A\times Y)\cup (X\times B)$. There
are three cases to consider. \begin{itemize}\item Suppose that
$c(s_*)\in A\times B$. Then $x(s_*)\in A$ and $y(s_*)\in B$. The
fact that $H$ and $J$ have the first coordinate independence
property with respect to $i$ and $j$ respectively means that each
coordinate of $Q$ is smooth, and so $Q$ is smooth. \item Suppose
that $c(s_*)\in A\times Y$, and that $y(s_*)\notin B$. Then at
each of the points $c(s_n)$, $(Q\circ c)(s_n)$ is given by the
second part of the definition of $Q$, for $n$ large enough. Since
$x(s_*)\in A$, the component of $Q$ involving $H$ is smooth, since
$H$ has the first coordinate independence property. For any $s$ in
a neighborhood of $s_*$,
$\alpha(\frac{\overline{u}(x(s))}{\overline{v}(y(s))})=0$. Thus,
the component of $Q$ involving $J$ is constant for $s$ in a
neighborhood of $s_*$, and so is smooth there. \item The case with
$c(s_*)\in X\times B$, and $x(s_*)\notin A$ is similar to the
second case above. \end{itemize}

For the last part of the theorem, suppose that $(u, H)$ represent
$(X, A)$ as an SDR pair. If we replace $u$ by $u'=\frac{1}{2}u$,
then $(u', H)$ also represent $(X, A)$ as an SDR pair. Making the
above constructions now with $u'$ in place of $u$, it follows that
$w(x, y)< 1$ for all $(x, y)$ and so $Q(1, x, y)\in (X\times
B)\cup (A\times Y)$. This completes the proof. \proofterminator

\section{Cofibrations} In this section, we show that for a
subspace $A\subseteq X$ that is closed in the underlying topology,
the inclusion $i: A\rightarrow X$ is a cofibration if and only if
$(X, A)$ is an SNDR pair.

\begin{mydf}\label{521} Let $i: A\rightarrow X$ be a cofibration. We call $i$
a {\bf cofibration with FCIP} if any homotopy extension can be
chosen to have the FCIP with respect to $i$. \end{mydf}Using the
equivalent formulation of the notion of cofibration, given by
Lemma \ref{463}, we may restate Definition \ref{521} as follows: A
cofibration $i: A\rightarrow X$ is a cofibration with the FCIP if
and only if the map $G$ that we may fill in to complete the
commutative diagram \[\xymatrix{(0\times X)\cup
(\mathbf{I}^-\times A)\ar[r]^{\hspace{13mm}h} \ar[d]_j & Y\\
\mathbf{I}^-\times X\ar@{-->}[ur]_G}\]may be chosen to have the
FCIP with respect to the inclusion $i$.

We have the following result, which corresponds to a similar
topological result.

\begin{lem}\label{522} A smooth map $i: A\rightarrow X$ is a cofibration
(with the FCIP) if and only if $(0\times X)\cup
(\mathbf{I}^-\times A)$ is a retract of $\mathbf{I}^-\times X$,
$($where the retraction $r: \mathbf{I}^-\times X\rightarrow
(0\times X)\cup (\mathbf{I}^-\times A)$ has the FCIP $)$.\end{lem}
\P In the one direction, suppose that $(0\times X)\cup
(\mathbf{I}^-\times A)$ is a retract of $\mathbf{I}^-\times X$. We
wish to complete the following diagram: \[\xymatrix{(0\times
X)\cup
(\mathbf{I}^-\times A)\ar[r]^{\hspace{13mm}h} \ar[d]_j & Y\\
\mathbf{I}^-\times X\ar@{-->}[ur]_G}.\]By hypothesis, there exists
$r: \mathbf{I}^-\times X\rightarrow (0\times X)\cup
(\mathbf{I}^-\times A)$ such that $r\circ j=1$. Define $G= h\circ
r$. If $r$ has the FCIP, then so does $h\circ r$.

Conversely, suppose that $i: A\rightarrow X$ is a cofibration
(with the FCIP). We may find a map $r$ such that the diagram
\[\xymatrix{(0\times X)\cup
(\mathbf{I}^-\times A)\ar[r]^{\hspace{13mm}1} \ar[d]_j & (0\times X)\cup (\mathbf{I}^-\times A)\\
\mathbf{I}^-\times X\ar@{-->}[ur]_r}\]commutes. Thus, $r\circ j=
1$. If $i$ is cofibration with the FCIP with respect to $i$, then
$r$ can be chosen to have the FCIP. \proofterminator

The next theorem shows the relationship between cofibrations,
retracts and SNDR pairs.

\begin{theo} Let $i: A\rightarrow X$ be an inclusion, with $A$
closed in the underlying topology of $X$. Then the following are
equivalent. \begin{enumerate}\item[{$(1)$}] The pair $(X, A)$ is
an SNDR pair. \item[{$(2)$}] There is a smooth retraction $r:
\mathbf{I}^-\times X\rightarrow (0\times X)\cup
(\mathbf{I}^-\times A)$ with the FCIP. \item[{$(3)$}] The map $i:
A\rightarrow X$ is a cofibration with the
FCIP.\end{enumerate}\end{theo} \P To show that $(1)$ and $(2)$ are
equivalent, note that the pair $(\mathbf{I}^-\times X, (0\times
X)\cup (\mathbf{I}^-\times A))$ is an SDR pair, as a consequence
of Lemma \ref{517} and Theorem \ref{519}. Let $(w, Q)$ be a
representation for the pair $(\mathbf{I}^-\times X, (0\times
X)\cup (\mathbf{I}^-\times A))$ as an SDR pair, and let $Q$ be
constructed as in Theorem \ref{519}. Define \[r:
\mathbf{I}^-\times X\rightarrow (0\times X)\cup
(\mathbf{I}^-\times A)\]by $r(t, x)= Q(1, t, x)$, where $(t, x)\in
\mathbf{I}^-\times X$. We observe that $r$ has the FCIP, since $Q$
has this property, and $Q$ has this property since each of its
components has this property.

The equivalence of $(2)$ and $(3)$ is Lemma \ref{522}.

We need only show that $(2)$ implies $(1)$. Let $r:
\mathbf{I}^-\times X\rightarrow (0\times X)\cup
(\mathbf{I}^-\times A)$ be a retraction with the FCIP with respect
to $i$. Define $H: \mathbf{I}\times X\rightarrow X$ by $H(t, x)=
(\pi_X\circ r)(\alpha(t), x)$, where $\pi_X$ is the projection
onto the second coordinate, and $\alpha: \mathbb{R}\rightarrow
\mathbb{R}$ is a braking function with the following properties:
$\alpha(t)= 0$ for $t\leq 0$, $\alpha(t)= 1$ for $t\geq
\frac{3}{4}$, and $0< \alpha(t) <1$ for $0< t< \frac{3}{4}$. This
braking function is necessary to ensure smoothness at the right
endpoint of the flattened unit interval $\mathbf{I}$. Smoothness
at the left endpoint is already taken care of by the fact that $r$
is defined in terms of the left flattened unit interval. The map
$H$ satisfies the following properties: \begin{itemize}\item $H$
has the FCIP since $r$ has this property. \item $H(0, x)=
(\pi_X\circ r)(0, x)= x$, for $x\in X$. \item $H(t, x)=
(\pi_X\circ r)(\alpha(t), x)= x$, for $x\in A$.\end{itemize}We now
construct $u: X\rightarrow \mathbf{I}$. Let $\pi_\mathbf{I}:
\mathbf{I}\times X\rightarrow \mathbf{I}$ denote the projection
onto $\mathbf{I}$. Define a smooth function $\beta:
\mathbb{R}\rightarrow \mathbb{R}$ by \[\beta(t)=
\left\{\begin{array}{ll} 0 & \mbox{if $t\leq 0$}\\
e^{-\frac{1}{t^2}} & \mbox{if $t>0$.}\end{array}\right.\]Now,
define $u: X\rightarrow \mathbf{I}$ by \[u(x)=
\frac{\int_0^1\beta(\alpha(t)-(\pi_\mathbf{I}\circ r)(1,
x)(\pi_\mathbf{I}\circ r)(\alpha(t),
x))dt}{\int_0^1\beta(\alpha(t))dt}.\]It is clear that $u$ is a
smooth mapping.

We now verify that $(u, H)$ represents $(X, A)$ as an SNDR pair.
\vspace{3mm}\newline $(1)$ Let $x\in A$. Clearly,
$(\pi_\mathbf{I}\circ r)(1, x)=1$ and $\pi_\mathbf{I}\circ
r)(\alpha(t), x)= \alpha(t)$, and so
$\int_0^1\beta(\alpha(t)-(\pi_\mathbf{I}\circ r)(1,
x)(\pi_\mathbf{I}\circ r)(\alpha(t), x))dt=0$. Thus, $u(x)= 0$,
for all $x\in A$. \vspace{3mm}\newline $(2)$ Suppose that $x\in
X-A$. Since $0\times (X-A)$ is open in the underlying topology on
$(0\times X)\cup (\mathbf{I}^-\times A)$, we may choose an open
neighborhood $W\subseteq 0\times (X-A)$ of $(0, x)$. Since $r$ is
continuous, there is a neighborhood $V\subseteq \mathbf{I}^-\times
X$ such that $r(V)\subseteq W\subseteq 0\times (X-A)$. Now,
consider the mapping $q_x: \mathbf{I}\rightarrow \mathbf{I}\times
X$, given by $q_x(t)= (\alpha(t), x)$, for each $x\in X$. This is
clearly smooth. Thus, there exists a neighborhood $U\subseteq
\mathbf{I}^-$ such that $q_x(U)\subseteq V$. In other words,
$U\times \{x\}\subseteq V$. So, we have $(\pi_\mathbf{I}\circ
r)(\alpha(t), x)= 0$, for all $t\in U$. Thus, we have \[u(x)=
\frac{\int_{I-U}\beta(\alpha(t)-(\pi_\mathbf{I}\circ r)(1,
x)(\pi_\mathbf{I}\circ r)(\alpha(t), x))dt+
\int_U\beta(\alpha(t))dt}{\int_0^1\beta(\alpha(t))dt}>0.\]Combining
this with part $(1)$, we deduce that $u^{-1}(0)= A$.
\vspace{3mm}\newline $(3)$ Suppose that $x$ is such that $u(x)<
1$. There must be a neighborhood $U$ of $\mathbf{I}$ such that
$(\pi_\mathbf{I}\circ r)(1, x)(\pi_\mathbf{I}\circ r)(\alpha(t),
x)>0$, for $t\in U$. Thus $(\pi_\mathbf{I}\circ r)(1, x)>0$, but
this implies that $r(1, x)\in \mathbf{I}\times A$, and hence $H(1,
x)\in A$. The proof is complete. \proofterminator

\section{The Mapping Cylinder} In this section we show that the
inclusion of $X$ into the flattened mapping cylinder $\mathbf{I}_f$
of a map $f: X\rightarrow Y$ is a cofibration with the FCIP.

\begin{theo}\label{531} Let $f: X\rightarrow Y$ be a smooth map. Then, the pair
$(\mathbf{I}_f, X)$ is an SNDR pair.\end{theo} \P Let $\alpha:
\mathbf{I}\rightarrow \mathbb{R}$ be a smooth braking function with
the following properties: $\alpha(t)=0$ if $0\leq t\leq
\frac{1}{4}$, $\alpha(t)= 1$ if $\frac{3}{4}\leq t\leq 1$, $0<
\alpha(t)< 1$, otherwise. Define two more braking functions
$\alpha_1, \alpha_2:\mathbf{I}\rightarrow \mathbb{R}$ as follows:
$\alpha_1(0)= 0$, $0< \alpha_1(t)< 1$ if $0< t< \frac{3}{4}$,
$\alpha_1(t)=1$ if $\frac{3}{4}\leq t\leq 1$, and $\alpha_2(t)= 0$
if $0\leq t\leq \frac{3}{4}$, $\alpha_2(t)=1$ if $\frac{7}{8}\leq
t\leq 1$. Now, define $u: \mathbf{I}_f\rightarrow \mathbf{I}$ by
$u([t, x])= \alpha_1(t)$ and $u([y])= 1$, for $(t, x)\in
\mathbf{I}\times X$ and $y\in Y$. Define $H: \mathbf{I}\times
\mathbf{I}_f\rightarrow \mathbf{I}_f$ by \[\left\{\begin{array}{ll}
H(s, [t, x])= [(1-\alpha(s))t+ \alpha(s)\alpha_2(t), x] & \mbox{if
$(t, x)\in \mathbf{I}\times X$}\\ H(s, [y])= [y] & \mbox{if $y\in
Y$}.\end{array}\right.\]

That $u$ is smooth comes from the fact that it is smooth when
restricted to each component of the coproduct $(\mathbf{I}\times
X)\sqcup Y$; it is thus smooth on the quotient $\mathbf{I}_f$.

To see that the map $H: \mathbf{I}\times \mathbf{I}_f\rightarrow
\mathbf{I}_f$ is smooth, note that since we are working in a
cartesian closed category, products commute with quotients, i.e. if
$q$ is quotient, then so is $1\times q$, where $1$ is an identity
map. Thus, we may think of $H$ as being defined on the space
\[\frac{(\mathbf{I}\times \mathbf{I}\times X)\sqcup (\mathbf{I}\times
Y)}{\sim}\]where $\sim$ is the identification $(t, 1, x)= (t, f(x))$
for $t\in \mathbf{I}$, and $x\in X$. Since $H$ is smooth when
restricted to each component of the coproduct $(\mathbf{I}\times
\mathbf{I}\times X)\sqcup (\mathbf{I}\times Y)$, $H$ is smooth on
the quotient $\mathbf{I}\times \mathbf{I}_f$.

We now verify that $(u, H)$ is a representation for $(\mathbf{I}_f,
X)$ as an SNDR pair.
\begin{itemize}\item $u^{-1}(0)= [0, x]=
i_0(X)$. \item $H(0, [t, x])= [t, x]$ and $H(0, [y])= [y]$. \item
$H(s, [0, x])= [0, x]$. \item If $u[t, x]< 1$, then $t< \frac{3}{4}$
and so $\alpha_2(t)=0$. Thus, $H(1, [t, x])= [0, x]$.\end{itemize}
This completes the proof. \proofterminator

Finally, we have the following important corollary.

\begin{cor} Given any smooth map $f: X\rightarrow Y$, the inclusion
$X\hookrightarrow \mathbf{I}_f$ is a cofibration with the FCIP.
\end{cor}

\section{The Exact Sequence of a Cofibration} Our aim in this
section is to show how one can use SNDR pairs to prove the existence
of the right exact Puppe sequence. We state the result in Theorem
\ref{541} and break the proof of the result up into a number of
lemmas. We follow the method used by Whitehead \cite{wh} for the
topological case.

Throughout this section we work in the category $\mathbb{FRL}_*$ of
pointed Fr|'olicher spaces, and basepoint preserving smooth maps.

\begin{theo}\label{541} Let $W$ be an object in $\mathbb{FRL}_*$, and suppose
that $i: A\hookrightarrow X$ is a cofibration in $\mathbb{FRL}_*$.
For any basepoint $x_0\in A\subseteq X$ there is a sequence
\[\xymatrix@1{{\ldots}\ar[r] & [\sum^{n+1}A, W]\ar[r]^{(\sum^nk)^*} &
[\sum^n\mathbf{T}_i, W]\ar[r]^{(\sum^nj)^*} & [\sum^nX,
W]\ar[r]^{(\sum^ni)^*} & [\sum^nA, W]\ar[r] & {\ldots}
\\{\ldots}\ar[r] & [\sum A, W]\ar[r]^{k^*} & [\mathbf{T}_i,
W]\ar[r]^{j^*} & [X, W]\ar[r]^{i^*} & [A, W]}\]which is an exact
sequence in $\mathbb{SETS}_*$, where $j: X\rightarrow \mathbf{T}_i$
is the inclusion discussed in Paragraphe \ref{45} and $k:
\mathbf{T}_i\rightarrow \sum A$ is the quotient map defined below.
\end{theo} It is, in fact, possible to prove that the sequence above
is an exact sequence of groups as far as $\sum A, W]$ and that the
morphisms to this point are group homomorphisms, but we shall not do
so here.

The $\mathbf{reduced (flattened ) suspension}$ of a pointed
Fr\"olicher space $X$ is defined as \[\sum X= (\mathbf{I}/\{0,
1\})\wedge X,\]where the reduced join is defined as for topological
spaces with the identified set taken as basepoint, and with $0$ the
basepoint of $\mathbf{I}$.

In this section, whenever we refer to the suspension of a space , we
mean the reduced flattened suspension defined above.

\begin{lem}\label{542} If $(x, A)$ is an SNDR pair and $p: X\rightarrow X/A$
the quotient map, then the sequence \[\xymatrix{A\ar[r]^i &
X\ar[r]^p & X/A}\]is right exact.\end{lem} \P To show that the given
sequence is right exact we must show that for any Fr\"olicher space
$W$ the following sequence is exact in $\mathbb{SETS}$:
\[\xymatrix{[X/A, W]\ar[r]^{p^*} & [X, W]\ar[r]^{i^*} & [A, W]}.\]It
is easy to see that im $p^*\subseteq \mbox{ker $i^*$}$. To see the
reverse inclusion, let $g: X\rightarrow W$ be an element of $[X,
W]$, with $g|_A\simeq w_0$ (rel $w_0$), where $w_0\in W$. Since
$\xymatrix{A\ar[r]^i & X}$ is an SNDR pair, the map $i$ is a
cofibration, and so we may extend $w_0$ to a smooth map $g':
X\rightarrow W$ such that $g'\simeq g$. But $g'$ is constant on $A$,
and so there exists a smooth map $g_1: X/A\rightarrow W$ such that
$p^*(g_1)= g'$. This shows that ker $i^*\subset \mbox{im $p^*$}$.
\proofterminator

\begin{lem}\label{543} For any smooth map $f: X\rightarrow Y$, the sequence
\[\xymatrix{X\ar[r]^f & Y\ar[r]^l & \mathbf{T}_f}\]is right exact,
where $l$ is the usual inclusion of $Y$ into the mapping cone; i.e.
$y\mapsto [y]\in \mathbf{T}_f$. \end{lem} \P One can show that there
is a homotopy commutative diagram \[\xymatrix{X\ar[r]^f\ar[dr]_i &
Y\ar[d]^j\ar[dr]^l\\ {} & \mathbf{I}_f\ar[r]_p &
\mathbf{T}_f}\]where $i, j$, and $l$ are the usual inclusions, and
$p$ is the quotient map that collapses away $\{0\}\times X$ to a
point. Since, by Theorem \ref{531}, $(\mathbf{I}_f, X)$ is an SNDR
pair, it follows from Lemma \ref{542} that the sequence
\[\xymatrix{X\ar[r]^i & \mathbf{I}_f\ar[r]^p & \mathbf{T}_f}\]is
right exact. It is fairly easy to show that $j: Y\rightarrow
\mathbf{I}_f$ is a homotopy equivalence. Therefore, the sequence
\[\xymatrix{X\ar[r]^f & Y \ar[r]^l & \mathbf{T}_f}\]is right exact.
\proofterminator

\begin{lem} For any smooth map $i: A\rightarrow X$, there is an
infinite right exact sequence \[\xymatrix{A\ar[r]^i & X\ar[r]^{i^1}
& \mathbf{T}_i\ar[r]^{i^2} & \ldots\ar[r]^{i^{n-1}} &
\mathbf{T}_{i^{n-2}}\ar[r]^{i^n} &
\mathbf{T}_{i^{n-1}}\ar[r]^{i^{n+1}} & \ldots}\]where $i^n$, $n\geq
1$, are inclusion maps. \end{lem} \P The pair $(\mathbf{T}_i, X)$ is
an SNDR pair. The representation for the pair $(\mathbf{I}_f, X)$ in
Theorem \ref{531} can be adapted to show this. One iterates the
procedure of Lemmas \ref{542} and \ref{543}. \proofterminator

One can easily see that there is an isomorphism between
$\mathbf{T}_i/X$ and $\sum A$. Define $q: \mathbf{T}_i\rightarrow
\sum A$ to be the map which identifies $X\subset \mathbf{T}_i$ to a
point, followed by the isomorphism $\mathbf{T}_i/X\rightarrow \sum
A$.

\begin{lem} The sequence \[X\xymatrix{X\ar[r]^{i^1} &
\mathbf{T}_i\ar[r]^q & \sum A}\] is right exact. \end{lem} \P As
noted above the pair $(\mathbf{T}_i, X)$ is an SNDR pair. We have
the commutative diagram \[\xymatrix{X\ar[r]^{i^1} &
\mathbf{T}_i\ar[r]^p\ar[dr]^q & \mathbf{T}_i/X\ar[d]^{q_0}\\ {} & {}
& \sum A}\]where $p: \mathbf{T}_i\rightarrow \mathbf{T}_i/X$ is the
identification map, and $q_0: \mathbf{T}_i/X\rightarrow \sum A$ is
an isomorphism. The top line of the diagram is right exact, by Lemma
\ref{542}, and so the sequence \[\xymatrix{X\ar[r]^{i^1} &
\mathbf{T}_i\ar[r]^q & \sum A}\]is right exact. \proofterminator

There is a commutative diagram \[\xymatrix{X\ar[r]^{i^1} &
\mathbf{T}_i\ar[r]^{i^2}\ar[dr]^q & \mathbf{T}_{i^1}\ar[d]^{q_1}\\
{} & {} & \sum A}\]where $q_1$ is a homotopy equivalence. ( See
Whitehead \cite{wh} for more details of this map. ) Using
commutative diagrams of this form, one can now proceed almost
exactly as one does in the topological situation, as in Whitehead
\cite{wh} for example, to get the following infinite right exact
sequence: \[\xymatrix{A\ar[r]^i & X\ar[r]^{i^1} &
\mathbf{T}_i\ar[r]^q & \sum A\ar[r]^{\sum i} & \sum X\ar[r]^{\sum
i^1} & \sum \mathbf{T}_i\ar[r]^{\sum q} & \ldots \\ {} & \ldots
\ar[r] & \sum^nA\ar[r]^{\sum^ni} & \sum^nX\ar[r]^{\sum^ni^1} &
\ldots}\]The definition of right exactness now gives us the exact
sequence of Theorem \ref{541}.

\end{document}